\documentclass[11pt]{article}
\usepackage{mathrsfs}
\usepackage{latexsym,lineno}
\usepackage{algorithm}
\usepackage{enumerate}
\usepackage{algorithmic}
\usepackage{epsfig}
\usepackage{color}
\usepackage{amsmath}\usepackage{fleqn}\usepackage{verbatim}
\usepackage{epsf}
\usepackage{amsthm}\usepackage{graphicx, float}\usepackage{graphicx}
\usepackage{amsfonts}
\usepackage{amssymb}\usepackage{graphpap}
\usepackage{epic}\usepackage{curves}
\usepackage{tikz}
\usepackage{subfigure}
\usepackage{mathrsfs}
\usepackage{pstricks}
\usepackage{hyperref}
\newtheorem{theorem}{Theorem}[section]

\newtheorem{lemma}{Lemma}[section]

\numberwithin{equation}{section}

\title{\bf Stability  of  generalized Tur\'an number \\for linear forests}
\author { Yisai Xue,\, Yichong Liu, \, Liying Kang\thanks{Corresponding author.   Email address: lykang@shu.edu.cn (L. Kang),  xys16720018@163.com (Y. Xue), lyc328az@163.com (Y. Liu). This work is  supported by
the National Nature Science Foundation of China (grant numbers 11871329, 11971298).}\\
{\small Department of Mathematics, Shanghai University,
Shanghai 200444, P.R. China}}

\date{}
\begin{document}

\maketitle

\begin{abstract}
Given a graph $T$ and a family of graphs $\mathcal{F}$, the generalized Tur\'an number of $\mathcal{F}$ is the maximum number of copies of $T$ in an $\mathcal{F}$-free graph on $n$ vertices, denoted by $ex(n,T,\mathcal{F})$. When $T = K_r$, $ex(n, K_r, \mathcal{F})$ is a function specifying the maximum possible number of $r$-cliques in an $\mathcal{F}$-free graph on $n$ vertices. A linear forest is a forest whose connected components are all paths and isolated vertices. Let $\mathcal{L}_{k}$ be the family of all linear forests of size $k$ without isolated vertices. In this paper, we obtained the maximum possible number of $r$-cliques in $G$, where $G$ is $\mathcal{L}_{k}$-free with minimum degree at least $d$. Furthermore, we give a stability version of the result. As an application of the stability version of the result, we obtain a clique version of the stability of the Erd\H{o}s-Gallai Theorem on matchings.

\bigskip

\medskip

\noindent{\bf AMS (2000) subject classification:}  05C35
\end{abstract}

\section{Introduction}

   Let $\mathcal{F}$ be a family of graphs. The \textit{Tur\'an number} of $\mathcal{F}$, denoted by $ex(n, \mathcal{F})$, is the maximum number of edges in a graph with $n$ vertices which does not contain any subgraph isomorphic to a graph in $\mathcal{F}$.
   When $\mathcal{F}=\{F\}$, we write $ex(n, F)$ instead of $ex(n, \{F\})$.
   The problem of determining Tur\'an number for assorted graphs traces its history back to 1907, when Mantel showed that $ex(n,K_3)=\lfloor\frac{n^2}{4}\rfloor$.
   In 1941, Tur\'an \cite{1941Turan} proved that if a graph does not contain a complete subgraph $K_r$, then the maximum number of edges it can contain is given by the Tur\'an-graph, a complete balanced $(r-1)$-partite graph.

   For a graph $G$ and $S,T\subseteq V(G)$, denote by $E_G(S,T)$ the set of edges between $S$ and $T$ in $G$, i.e., $E_G(S,T)=\{uv\in E(G)\colon\, u\in S, v\in T\}$.
   Let $e_G(S,T)=|E_G(S,T)|$.
   If $S=T$, we use $e_G(S)$ instead of $e_G(S,S)$.
   For  a vertex $v\in V(G)$,  the {\it degree} of  $v$, written as $d_G(v)$ or simply $d(v)$,
is the number of edges incident with $v$. We use $d_T(v)$ instead of $e_G(S,T)$ when $S=\{v\}$.
   For any $U \subseteq V(G)$,  let $G[U]$ be the subgraph   induced by $U$ whose edges are precisely the edges of $G$ with
both ends in $U$.

 Let $G$ be a graph of order $n$, $P$ a property defined on $G$, and $k$ a positive integer.
   A property $P$ is said to be \textit{$k$-stable}, if whenever $G+uv$ has the property $P$ and $d_G(u) + d_G(v) \geq k$, then $G$ itself has the property $P$.
The $k$-\textit{closure} of a graph $G$ is the (unique) smallest graph $G'$ of order $n$ such that $E(G) \subseteq E(G')$ and $d_{G'}(u)+d_{G'}(v)<k$ for all $u v \notin E(G')$.
   The $k$-closure can be obtained from $G$ by a recursive procedure of joining nonadjacent vertices with degree-sum at least $k$. In particular, if $G'=G$, we say that $G$ is \textit{stable under taking} $k$-closure.
   Thus, if $P$ is $k$-stable and the $k$-closure of $G$ has property $P$, then $G$ itself has property $P$.

   For a natural number $\alpha$ and a graph $G$, the $\alpha$-\textit{disintegration} of a graph $G$ is the process of iteratively removing from $G$ the vertices with degree at most $\alpha$ until the resulting graph has minimum degree at least $\alpha+1$ or is empty. The resulting subgraph $H=H(G, \alpha)$ will be called the $(\alpha+1)$-\textit{core} of $G$.
It is well known that $H(G, \alpha)$ is unique and does not depend on the order of vertex deletion (for instance, see \cite{1996P}).
   The \textit{matching number} $\nu(G)$ is the number of edges in a maximum matching of $G$.

The $n$-vertex graph $H(n, k, a)$ is defined as follows.
   The vertex set of $H(n, k, a)$ is partitioned into three sets $A, B, C$ such that $|A|=a,|B|=k-2a,|C|=n-k+a$, and the edge set of $H(n, k, a)$ consists of all edges between $A$ and $C$ together with all edges in $A \cup B$.
   Let $H^+(n, k, a)$ and $H^{++}(n, k, a)$ be the graph obtained by adding one edge and two independent edges in $C$ of $H(n, k, a)$, respectively.
   The number of $r$-cliques in $H(n, k, a)$ is denoted by $h_{r}(n, k, a):=\binom{k-a}{r}+(n-k+a)\binom{a}{r-1}$, where $h_{r}(n, k, 0)=\binom{k}{r}$.
%   Note that $h_{2}(n, k, a)$ is the number of edges of $H(n, k, a)$, and we write $h(n, k, a)$  instead of $h_{2}(n, k, a)$.
%   Many Tur\'an-type problems have an extremal graph that is isomorphic to $H(n, k, a)$.
%   For example, the extremal graph in Theorem \ref{1959E} is either $H(n, 2k-1, 0)$ or $H(n, 2k-1, k-1)$.

    A \textit{linear forest} is a forest whose connected components are all paths and isolated vertices.
   %This type of linear forests is also well studied (see Hu et al. \cite{2001Hu}).
   Let  $\mathcal{L}_{k}$ be the family of all linear forests of size $k$ without isolated vertices.
   In \cite{2019Wang}, Wang and Yang  proved that $ex\left(n ; \mathcal{L}_{n-k}\right)=\binom{n-k}{2}+O\left(k^{2}\right)$ when $n\geq 3k$.
   Later, Ning and Wang \cite{2020Ning} completely determined the Tur\'an number $ex\left(n ; \mathcal{L}_{k}\right)$ for all $n>k$.

\begin{figure}
\begin{center}
\begin{tikzpicture}
\filldraw [color=black, fill=blue!15, very thick] (0,0) circle (50pt)node [left=2pt]{\Large $B$};
\filldraw [color=black, fill=blue!15, very thick] (1,0) circle (20pt)node {\Large $A$};
\draw[color=black, fill=white, , very thick] (4.5,0) ellipse(0.8 and 1.6)node {\Large $C$};
\draw  [black,very thick](1,0.7)--(4.5,1.6);
\draw  [black,very thick](1,-0.7)--(4.5,-1.6);
\end{tikzpicture}
\end{center}\caption{$H(n,k,a)$}
\end{figure}

\begin{theorem}[Ning and Wang \cite{2020Ning}]\label{span}
   For any integers $n$ and $k$ with $1 \leq k \leq n-1$, we have
$$
ex\left(n,\mathcal{L}_{k}\right)=\max \left\{h_2\left(n, k, 0\right),h_2\left(n, k, \left\lfloor\frac{k-1}{2}\right\rfloor\right)\right\}.
$$
\end{theorem}

   Given a graph $T$ and a family of graphs $\mathcal{F}$, the \textit{generalized Tur\'an number} of $\mathcal{F}$ is the maximum number of copies of $T$ in an $\mathcal{F}$-free graph on $n$ vertices, denoted by $ex(n,T,\mathcal{F})$.
   Note that $ex(n, K_2, \mathcal{F})=ex(n, \mathcal{F})$.
   The problem to estimate generalized Tur\'an number has
received a lot of attention.
   In 1962, Erd\H{o}s \cite{1962E} generalized the classical result of Tur\'an by determining the exact value of $ex(n,K_r,K_t)$.
   Luo \cite{2018Luo} determined the upper bounds on $ex(n,K_r,P_{k})$ and $ex(n,K_r,\mathcal{C}_{\geq k})$, where $\mathcal{C}_{\geq k}$ is the family of all cycles with length at least $k$.
   In \cite{2020Gerbner}, Gerbner, Methuku and Vizer  investigated the function $ex(n,T,kF)$, where $kF$ denotes $k$ vertex disjoint copies of a fixed graph $F$.
   The systematic study of $ex(n,T,\mathcal{F})$ was initiated by Alon and Shikhelman \cite{2016Alon}.
   Recently, Zhang, Wang and Zhou \cite{2021Zhang} determined the exact values of $ex(n,K_r,\mathcal{L}_{k})$ by using  the shifting method.

\begin{theorem}[Zhang, Wang and Zhou \cite{2021Zhang}]\label{2021Zhang}
   For any $r\geq 2$ and $n\geq k+1$,
$$ex(n,K_r,\mathcal{L}_{k})=\max \left\{h_r\left(n, k, 0\right), h_r\left(n, k, \left\lfloor\frac{k-1}{2}\right\rfloor\right)\right\}.$$
\end{theorem}

   Let $N_r(G)$ denote the number of $r$-cliques in $G$.
   When $T = K_r$, $ex(n, K_r, \mathcal{F})$ is a function specifying the maximum possible number of $r$-cliques in an $\mathcal{F}$-free graph on $n$ vertices.
   We extend Theorem \ref{2021Zhang} as follows.

\begin{theorem}\label{clique}
   Let $G$ be an $\mathcal{L}_{k}$-free graph on $n$ vertices with minimum degree $d$ and $d\leq \lfloor\frac{k-1}{2}\rfloor$.
   Then $$N_r(G)\leq \max \left\{h_r\left(n, k, d\right), h_r\left(n, k, \left\lfloor\frac{k-1}{2}\right\rfloor\right)\right\}.$$
   The graphs $H(n,k,d)$ and $H\big(n,k,\lfloor\frac{k-1}{2}\rfloor\big)$ show that this bound is sharp.
\end{theorem}

   Many extremal problems have the property that there is a unique extremal example, and moreover any construction of close to maximum size is structurally close to this extremal example.
  In \cite{2019F}, F\"uredi, Kostochka, and Luo studied the maximum number of cliques in non-$\ell$-hamiltonian graphs,  where the property non-$\ell$-hamiltonian is $(n+\ell)$-stable.
   Actually, they not only asked to determine the maximum number of cliques in graphs having a stable property $P$, but also asked to prove a stability version of it.
   Motivated by the question  proposed by F\"uredi, Kostochka, and Luo \cite{2019F}, we give the following result  which is the stability version of Theorem \ref{clique}.

%\begin{theorem}\label{l_stability}
%   Let $G$ be an $\mathcal{L}_{k}$-free graph on $n$ vertices with minimum degree at least $d$.
%   If $n> k^3$, and $$N_r(G)>\max \left\{h_r(n, k, d), h_r\left(n, k, \left\lfloor\frac{k-3}{2}\right\rfloor\right)\right\},$$ then $G$ is a subgraph of $H(n, k, \lfloor\frac{k-1}{2}\rfloor)$.
%\end{theorem}

\begin{theorem}\label{stab2}
   Let $G$ be an $\mathcal{L}_{k}$-free graph on $n$ vertices with minimum degree at least $d$.
   If $n > k^5$, $r\leq \lfloor\frac{k-3}{2}\rfloor$  and $$N_r(G)>\max \left\{h_r(n, k, d),h_r\left(n, k, \left\lfloor\frac{k-5}{2}\right\rfloor\right)\right\},$$ then \\
   (i) $G$ is a subgraph of the graph $H\left(n, k, \left\lfloor\frac{k-1}{2}\right\rfloor\right)$, $H\left(n, k, \left\lfloor\frac{k-3}{2}\right\rfloor\right)$ or $H^+\left(n, k-1, \left\lfloor\frac{k-3}{2}\right\rfloor\right)$ if $k$ is odd;\\
   (ii) $G$ is a subgraph of the graph $H\left(n, k, \left\lfloor\frac{k-1}{2}\right\rfloor\right)$, $H\left(n, k, \left\lfloor\frac{k-3}{2}\right\rfloor\right)$, $H^+\left(n, k-1, \left\lfloor\frac{k-3}{2}\right\rfloor\right)$ or $H^{++}\left(n, k-2, \left\lfloor\frac{k-3}{2}\right\rfloor\right)$ if $k$ is even.
\end{theorem}

%\begin{theorem}\label{l_stability2}
%   Let $G$ be an $\mathcal{L}_{k}$-free graph on $n$ vertices with minimum degree at least $d$.
%   If $n \geq k+d+3$ and $$N_r(G)>\max \left\{h_r(n, k, d+1), h_r\left(n, k, \left\lfloor\frac{k-1}{2}\right\rfloor\right)\right\},$$ then $G$ is a subgraph of $H(n, k, d)$.
%\end{theorem}

   In 1959, Erd\H{o}s and Gallai \cite{1959E} determined  the maximum numbers of edges in an $n$-vertex graph with $\nu(G)\leq k$.

\begin{theorem}[Erd\H{o}s-Gallai Theorem \cite{1959E}]\label{1959E}
    Let $G$ be a graph on $n$ vertices.
    If $\nu(G) \leq k$, then
$$e(G)\leq \max \left\{h_{2}(n, 2k+1, 0),h_{2}(n, 2k+1, k)\right\}.$$
\end{theorem}

   In \cite{2020Duan}, Duan et al. extended Erd\H{o}s-Gallai Theorem as follows.

\begin{theorem}[Duan et al. \cite{2020Duan}]\label{clique11}
   If $G$ is a graph with $n\geq 2k+2$ vertices, minimum degree $d$, and $\nu(G) \leq k$,
   then $$N_r(G)\leq\max \left\{h_r\left(n, 2k+1, d\right), h_r\left(n, 2k+1, k\right)\right\}.$$
\end{theorem}

   As an application of our result, we give the stability version of Theorem \ref{clique11} for $2\leq r\leq k-1$.

\begin{theorem}\label{thm11}
   Let $G$ be a graph on $n$ vertices with $\delta(G) \geq d$ and $\nu(G) \leq k$.
   If $r\leq k-1$,   $n > (2k+1)^5$ and
   $$N_r(G)>\max \left\{h_{r}(n, 2k+1, d),h_{r}(n, 2k+1, k-2)\right\},$$
   then $G$ is a subgraph of $H(n, 2k+1, k)$ or $H(n, 2k+1, k-1)$.
\end{theorem}

\section{The maximum number of cliques in $\mathcal{L}_{k}$-free graphs with given minimum degree}

   The closure technique, which is initiated by Bondy and Chv\'atal \cite{1976Bondy} in 1976, played a crucial role in the proof of Theorem \ref{clique}.
%   Let $G$ be a graph of order $n$, $P$ a property defined on $G$, and $k$ a positive integer.
%   A property $P$ is said to be $k$-stable, if whenever $G+uv$ has the property $P$ and $d_G(u) + d_G(v) \geq k$, then $G$ itself has the property $P$.
In \cite{2020Ning}, Ning and Wang proved the property $\mathcal{L}_{k}$-free is $k$-stable.

\begin{lemma}[\cite{2020Ning}]\label{closure}
   Let $G$ be a graph on $n$ vertices. Suppose that $u, v \in V(G)$ with $d(u)+d(v) \geq k$. Then $G$ is $\mathcal{L}_{k}$-free if and only if $G+u v$ is $\mathcal{L}_{k}$-free.
\end{lemma}

\noindent\textbf{Proof of Theorem \ref{clique}.}
   Suppose, by way of contradiction, that $G$ is an $\mathcal{L}_{k}$-free graph with $N_{r}(G)>\max \left\{h_r\left(n, k, d\right), h_r\left(n, k, \left\lfloor\frac{k-1}{2}\right\rfloor\right)\right\}$.
   Let $G^{\prime}$ be the $k$-closure of $G$.
   Then  Lemma \ref{closure} implies that $G'$ is $\mathcal{L}_{k}$-free.
   Obviously, $\delta\left(G^{\prime}\right) \geq \delta(G)=d$.
   Let $H_{1}$ denote the $\lfloor\frac{k+1}{2}\rfloor$-core of $G'$, i.e., the resulting graph of applying $\lfloor\frac{k-1}{2}\rfloor$-disintegration to $G'$.

\noindent\textbf{Claim 1.} $H_{1}$ is nonempty.

\noindent\textbf{Proof.} Suppose $H_{1}$ is empty.
   Since one vertex is deleted at each step during the process of $\lfloor\frac{k-1}{2}\rfloor$-disintegration,  that destroys at most $\binom{\lfloor\frac{k-1}{2}\rfloor}{r-1}$ cliques of size $r$.
   The number of $K_{r}$'s contained in the last $\lceil\frac{k+1}{2}\rceil$ vertices is at most $\binom{\lceil\frac{k+1}{2}\rceil}{r}$.
   Therefore,
$$
\begin{aligned}
N_{r}\left(G^{\prime}\right) & \leq \binom{\left\lceil\frac{k+1}{2}\right\rceil}{r} + \left(n-\left\lceil\frac{k+1}{2}\right\rceil\right)\binom{\lfloor\frac{k-1}{2}\rfloor}{r-1}\\
 & = h_{r}\left(n, k, \left\lfloor\frac{k-1}{2}\right\rfloor\right) \\
 & \leq \max \left\{h_r\left(n, k, d\right), h_r\left(n, k, \left\lfloor\frac{k-1}{2}\right\rfloor\right)\right\},
\end{aligned}
$$
contradicting to the assumption of $G^{\prime}$, the claim follows.\qed

\noindent\textbf{Claim 2.} $H_{1}$ is a clique.

\noindent\textbf{Proof.}
Note that $d_{G^{\prime}}(u)\geq \lfloor\frac{k+1}{2}\rfloor$ for any vertex $u$ in $H_1$.
   Since $G^{\prime}$ is closed under taking $k$-closure, $H_1$ is a clique.\qed

   Let $t=\left|V\left(H_{1}\right)\right|$.
   Now we estimate the range of $t$.

\noindent\textbf{Claim 3.} $\lfloor\frac{k+3}{2}\rfloor\leq t \leq k-d$.

\noindent\textbf{Proof.}
   As $H_{1}$ is a clique and $d_{H_{1}}(u) \geq \lfloor\frac{k+1}{2}\rfloor$ for any vertex $u$ in $H_{1}$, we get $t \geq \lfloor\frac{k+3}{2}\rfloor$.
   If $t \geq k-d+1$, then $d_{G^{\prime}}(u) \geq d_{H_{1}}(u) = t-1\geq k-d$ for any vertex $u$ in $H_1$.
   Let $v$ be any vertex in $V\left(G^{\prime}\right) \backslash V\left(H_{1}\right)$.
   Notice that $d_{G^{\prime}}(v) \geq d_{G}(v) \geq d$ and $d_{G^{\prime}}(u)+d_{G^{\prime}}(v) \geq k-d+d=k$.
   Since $G^{\prime}$ is the $k$-closure of $G$, $v$ is adjacent to $u$.
   Then $G^{\prime}$ contains a $P_{k+1}$, which is  a contradiction.
   Thus $\lfloor\frac{k+3}{2}\rfloor \leq t \leq k-d$.\qed

   Let $H_{2}$ be the $(k+1-t)$-core of $G^{\prime}$.
   Since $t \geq \lfloor\frac{k+3}{2}\rfloor$, we obtain $k+1-t \leq \lfloor\frac{k+1}{2}\rfloor$.
   Therefore, $H_{1} \subseteq H_{2}$.

\noindent\textbf{Claim 4.} $H_{1} \neq H_{2}$.

\noindent\textbf{Proof.}
   Suppose $H_{1}=H_{2}$.
   Then $\left|V\left(H_{2}\right)\right|=t$.
   Since each step during the process of $(k-t)$-disintegration  destroys at most $\binom{k-t}{r-1}$ cliques of size $r$,
   we have $N_{r}\left(G^{\prime}\right) \leq\binom{t}{r}+(n-t)\binom{k-t}{r-1}=h_{r}(n, k, k-t)$.
   Note that $d \leq k-t\leq \lceil\frac{k-3}{2}\rceil$ from Claim 3.
   By the convexity of $h_{r}(n, k, k-t)$, we have $N_{r}\left(G^{\prime}\right) \leq \max \left\{h_{r}(n, k, d), h_{r}(n, k, \lceil\frac{k-3}{2}\rceil)\right\}\leq \max \left\{h_r\left(n, k, d\right), h_r\left(n, k, \left\lfloor\frac{k-1}{2}\right\rfloor\right)\right\}$, a contradiction.
   Thus the claim follows.\qed

   By Claim 4,  $H_{1}$ is a proper subgraph of $H_{2}$.
   This implies that there are non-adjacent vertices $u$ and $v$ such that $u \in V(H_{1})$ and $v \in V(H_{2}) \backslash V(H_{1})$.
   We have $d_{G^{\prime}}(u)+d_{G^{\prime}}(v) \geq t-1+(k+1-t)=k$.
   As $G^{\prime}$ is stable under taking $k$-closure, $u$ must be adjacent to $v$.
   We obtained a contradiction.

   It is easy to see that graphs $H(n,k,d)$ and $H(n,k,\lfloor\frac{k-1}{2}\rfloor)$ are $\mathcal{L}_{k}$-free. Then either $H(n,k,d)$  or $H(n,k,\lfloor\frac{k-1}{2}\rfloor)$ obtains the bound.
   The theorem is proved.\qed

\section{Stability on $\mathcal{L}_{k}$-free graphs}

\subsection{Proof of Theorem \ref{stab2}}

   Let $G$ be a graph on $n$ vertices.
   If there are at least $s$ vertices in $V(G)$ with degree at most $q$, then we say $G$ has $(s, q)$-\textit{P\'osa property}.
   If $G$ has $(s, q)$-P\'osa property and $n \geq s+q$, then we can check that
   $$N_r(G) \leq\binom{n-s}{r}+s \binom{q}{r-1}.$$

%   The following lemma is in the same spirit of Lemma 11 in \cite{2019F}, which is a key ingredient in our proof of the stability theorems.
   The following two lemmas show the relationship between the $k$-stable property and the P\'osa property.
   With the help of these two lemmas, we can approximate the structure of $k$-closure of a graph.

\begin{lemma}\label{posa}
   Let $n \geq k+1$.
   Assume property $P$ is $k$-stable and the complete graph $K_n$ has the property $P$.
   Suppose $G$ is a graph on $n$ vertices with minimum degree at least $d$.
   If $G$ does not have property $P$, then there exists an integer $q$ with $d \leq q \leq \frac{k-1}{2}$ such that G has $(n-k+q, q)$-P\'osa property.
\end{lemma}

\noindent\textbf{Proof.}
   Let $G^{\prime}$ be the $k$-closure of $G$ and $d_{G^{\prime}}\left(v_{1}\right), d_{G^{\prime}}\left(v_{2}\right), \cdots, d_{G^{\prime}}\left(v_{n}\right)$ be the degree sequence of $G^{\prime}$ such that $d_{G^{\prime}}\left(v_{1}\right)$ $\geq d_{G^{\prime}}\left(v_{2}\right) \geq \cdots \geq d_{G^{\prime}}\left(v_{n}\right)$.
   Clearly, $G'$ is not a complete graph.
   Otherwise $G'$ has property $P$, so does $G$, a contradiction.

   Let $v_{i}$ and $v_{j}$ be two non-adjacent vertices in $G^{\prime}$ with $1 \leq i<j \leq n$ and $d_{G^{\prime}}\left(v_{i}\right)+d_{G^{\prime}}\left(v_{j}\right)$ as large as possible.
   Obviously, $d_{G^{\prime}}\left(v_{i}\right)+d_{G^{\prime}}\left(v_{j}\right) \leq k-1$.
   Let $S$ be the set of vertices in $V \backslash\{v_i\}$ which are not adjacent to $v_{i}$ in $G'$.
   By the choice of $v_j$, we have $d_{G^{\prime}}(v) \leq d_{G^{\prime}}\left(v_{j}\right)$ for any $v \in S$.
   Then
$$
|S|=n-1-d_{G^{\prime}}\left(v_{i}\right) \geq n-k+d_{G^{\prime}}\left(v_{j}\right) .
$$
   There are at least $n-k+d_{G^{\prime}}\left(v_{j}\right)$ vertices in $V\left(G^{\prime}\right)$ with degree at most $d_{G^{\prime}}\left(v_{j}\right)$.
   Let $q=d_{G^{\prime}}\left(v_{j}\right)$.
   Then $G'$ has $(n-k+q, q)$-P\'osa property.
   Moreover, since $d_{G^{\prime}}\left(v_{i}\right) \geq d_{G^{\prime}}\left(v_{j}\right)$ and $d_{G^{\prime}}\left(v_{i}\right)+d_{G^{\prime}}\left(v_{j}\right) \leq k-1$, it follows that $q=d_{G^{\prime}}\left(v_{j}\right) \leq \frac{k-1}{2}$.
   Since $G$ is a subgraph of $G^{\prime}$ and
$$
d_{G^{\prime}}\left(v_{j}\right) \geq \delta\left(G^{\prime}\right) \geq \delta(G) \geq d,
$$
   we complete the proof.\qed

   The following lemma gives a structural characterization of graphs with P\'osa property.

\begin{lemma}\label{bipartite}
   Suppose G has $n$ vertices and is stable under taking $k$-closure.
   Let $q$ be the minimum integer such that $G$ has $(n-k+q, q)$-P\'osa property and $q \leq \frac{k-1}{2}$.
   If $T$ is the set of vertices in $V(G)$ with degree at least $k-q$ and $T^{\prime}=V(G) \backslash T$, then $G\left[T, T^{\prime}\right]$ is a complete bipartite graph.
\end{lemma}

\noindent\textbf{Proof.}
   Assume that $G\left[T, T^{\prime}\right]$ is not a complete bipartite graph.
   Choose two non-adjacent vertices $u \in T$ and $v \in T^{\prime}$ such that $d(u)+d(v)$ is as large as possible.
   Clearly, $d(u)+d(v) \leq k-1$ and $T$ forms a clique in $G$ as $G$ is stable under taking $k$-closure.
   Now denote by $S$ the set of vertices in $V \backslash\{u\}$ which are not adjacent to $u$ in G.
Clearly, for any $v^{\prime} \in S, d\left(v^{\prime}\right) \leq d(v)$ and
$$
|S|=n-1-d(u) \geq n-k+d(v).
$$
   Since $d(u) \geq k-q$ and $d(u)+d(v) \leq k-1, d(v) \leq q-1$.
   Let $q^{\prime}=d(v) \leq q-1$.
   We have at least $n-k+q^{\prime}$ vertices in $V(G)$ with degree at most $q^{\prime}$.
   Then $G$ has $(n-k+q', q')$-P\'osa property with $q'<q$, which contradicts the minimality of $q$.
   The lemma follows.\qed

%\begin{figure}\label{ct}
%\begin{center}
%\begin{tikzpicture}
%\filldraw [black] (0,0) circle (1.5pt)node [left=1 pt]{$v_0$};
%\filldraw [black] (1,0) circle (1.5pt);
%\filldraw [black] (2,0) circle (1.5pt);
%\filldraw [black] (4,0) circle (1.5pt);
%\filldraw [black] (5,0) circle (1.5pt)node [right=1 pt]{$v_t$};
%\filldraw [black] (3-0.3,0) circle (0.5pt);
%\filldraw [black] (3,0) circle (0.5pt);
%\filldraw [black] (3+0.3,0) circle (0.5pt);
%\filldraw [black] (0.7,1) circle (1.5pt);
%\filldraw [black] (0.9,1) circle (0.5pt)node [above=1 pt]{$\Delta-2$};
%\filldraw [black] (1.1,1) circle (0.5pt);
%\filldraw [black] (1.3,1) circle (1.5pt);
%\filldraw [black] (1.7,1) circle (1.5pt);
%\filldraw [black] (1.9,1) circle (0.5pt)node [above=1 pt]{~~$\Delta-2$};
%\filldraw [black] (2.1,1) circle (0.5pt);
%\filldraw [black] (2.3,1) circle (1.5pt);
%\filldraw [black] (3.7,1) circle (1.5pt);
%\filldraw [black] (3.9,1) circle (0.5pt)node [above=1 pt]{~$\Delta-2$};
%\filldraw [black] (4.1,1) circle (0.5pt);
%\filldraw [black] (4.3,1) circle (1.5pt);
%\draw [black](0,0)--(1,0);
%\draw [black](1,0)--(2,0);
%\draw [black](4,0)--(5,0);
%\draw [black](1,0)--(0.7,1);
%\draw [black](1,0)--(1.3,1);
%\draw [black](2,0)--(1.7,1);
%\draw [black](2,0)--(2.3,1);
%\draw [black](4,0)--(3.7,1);
%\draw [black](4,0)--(4.3,1);
%\draw [black](0,0) .. controls (2,-1) and  (3,-1) .. (5,0);
%%\draw  [black] (2.5, 0) ellipse (2.5 and 0.5);
%%\draw[black] (5,0) arc (0:-180:2.5);
%\end{tikzpicture}
%\end{center}\caption{$C(t,\Delta)$}
%\end{figure}

   Let $g(k,\Delta)$ be the maximum number of edges in a graph such that the size of linear forests is at most $k$ and the maximum degree is at most $\Delta$.
   The following lemma estimates the upper bound of $g(k,\Delta)$.

\begin{lemma}\label{k2}
For $k\geq 1$ and $\Delta\geq 3$, \\
(i) ~$g(k,2)\leq \frac{3}{2}k.$\\
(ii) $g(k,\Delta)\leq k(\Delta-1).$
\end{lemma}

\noindent\textbf{Proof of (i).}
   Let $G$ be an $\mathcal{L}_{k+1}$-free graph with $e(G)=g(k,2)$ and $\Delta(G)\leq 2$.
   Clearly, $g(1,2)=1$ and $g(2,2)=3$.
   Now  suppose that $k\geq 3$.
   Since the maximum degree is at most 2, each nontrivial component is either a path or a cycle.
   We claim that each component with at least 3 vertices is a cycle.
   If not, we add an edge between the two ends of the path and the resulting graph is still $\mathcal{L}_{k+1}$-free, which contradicts the maximality of $G$.
   If there is a component consisting of exactly one edge, we replace this edge and a component $C_{\ell}$ in $G$ with $C_{\ell+1}$.
   Then the resulting graph is still $\mathcal{L}_{k+1}$-free and the number of edges is equal to $g(k,2)$.
   Therefore, we can further assume that each nontrivial component of $G$ is a cycle.

   Let $C_{k_1}$, \ldots, $C_{k_t}$ be the nontrivial components of $G$.
   Then $k=(k_1-1) + \cdots + (k_t-1)$ and $e(G)=k_1+\cdots+k_t=k+t$.
   Note that $k_i-1\geq 2$, we have $t\leq \frac{k}{2}$.
   Thus $g(k,2)=e(G)\leq \frac{3}{2}k$.\qed

\noindent\textbf{Proof of (ii).}
   We use induction on $k$.
   It is easy to check that $g(1,\Delta)=1$ and $g(2,\Delta)=\Delta$.
   Thus lemma holds for $k=1,2$.
   Suppose that the lemma holds for all $k'<k$.
   Let $G$ be an $\mathcal{L}_{k+1}$-free graph with  $\Delta(G)\leq \Delta$.
   Let $P=v_0v_1\cdots v_{t}$ be the longest path in $G$ and $B=V(G)\backslash V(P)$.
   Then $G[B]$ is $\mathcal{L}_{k+1-t}$-free and $e(G[B])\leq (k-t)(\Delta-1)$ by the induction hypothesis.

   	Since $P$ is the longest path in $G$, $d_B(v_0)=d_B(v_t)=0$ and $d_B(v_i)\leq \Delta-2$ for $1\leq i\leq t-1$.
   	Thus,
$$
\begin{aligned}
  e(G[V(P)])+e_G[V(P), B])
  & = \frac{1}{2}\left(\sum\limits_{i=0}^t d_G(v_i)+\sum\limits_{i=0}^t d_B(v_i)\right)\\
  & \leq \frac{1}{2}\left((t+1)\Delta+(t-1)(\Delta-2)\right)\\
  & = t(\Delta-1)+1
\end{aligned}
$$
   The equality holds only if   $d_G(v_0)=\cdots=d_G(v_t)=\Delta$, $d_B(v_1)=\cdots=d_B(v_{t-1})=\Delta-2$ and $d_B(v_0)=d_B(v_t)=0$  hold simultaneously, which is impossible.
   Therefore, $e(G[V(P)])+e_G[V(P), B])\leq t(\Delta-1)$.
   Moreover, we have
   $$
\begin{aligned}
   e(G)&=e(G[B])+e(G[V(P)])+e_G[V(P), B])\\
   &\leq (k-t)(\Delta-1)+t(\Delta-1)\\
   &\leq k(\Delta-1).
\end{aligned}
$$
\qed

%   Let $T'$ be the graph after using operation on $T$.
%   Notice that $P$ is the longest path in $G$, there is no vertex in $G[B]$ is adjacent to $v_0$ and $v_t$.
%   Since the maximum degree of $G$ is at most $\Delta$, there are at most $\Delta-2$ vertices in $V(G)\backslash\{v_{i-1},v_{i+1}\}$ are adjacent to $v_i$, where $1\leq i\leq t-1$.
%   Let $C(t,\Delta)$ denote the graph obtained by attaching $t-2$ pendent edges to each of the consecutive $t-1$ vertices on $C_{t+1}$ (see Figure 2).
%   Then $T'$ is a subgraph of $C(t,\Delta)$ and $e(T)=e(T')\leq e(C(t,\Delta))\leq t(\Delta-1)$.
%   Therefore, $e(G)=e(G[B])+e(T)\leq (k-t)(\Delta-1)+t(\Delta-1)=k(\Delta-1)$.\qed

\noindent\textbf{Remark.}
   The graph consisting of $k/3$ pairwise disjoint $K_4$'s shows the bound in Lemma \ref{k2} (ii) is sharp when  $3$ divides $k$ and $\Delta=3$.

For integers $m, l, r$, the following combinatorial identity is well-known.
\begin{eqnarray}\label{3.1}
\binom{m+l}{r}=\sum_{j=0}^{r}\binom{m}{j}\binom{l}{r-j}
\end{eqnarray}

   The following lemma bounds the number of $r$-cliques by the number of edges.

\begin{lemma}[\cite{2021Chakraborti}]\label{cliques}
Let $r \geq 3$ be an integer, and let $x \geq r$ be a real number. Then, every graph with exactly $\binom{x}{2}$ edges contains at most $\binom{x}{r}$ cliques of order $r$.
\end{lemma}

   For two disjoint vertex sets $T$ and $T'$ of $G$,  we use $N_r^i\left(T, T'\right)$ and $N_r^{\geq i}\left(T, T'\right)$ to denote
    the number of $r$-cliques in $G[T,T']$ that contain exactly $i$ vertices and at least $i$ vertices in $T'$,  respectively.

\noindent\textbf{Proof of Theorem \ref{stab2}.}
   Let $G^{\prime}$ be the $k$-closure of $G$.
   Then $G'$ is $\mathcal{L}_{k}$-free from Lemma \ref{closure}.
   By Lemma \ref{posa}, there exists an integer $q$ with $d \leq q \leq \lfloor\frac{k-1}{2}\rfloor$ such that $G^{\prime}$ has $(n-k+q, q)$-P\'osa property.
   Furthermore, we assume $q$ is as small as possible.
   Then either $q=\lfloor\frac{k-1}{2}\rfloor$ or $q=\lfloor\frac{k-3}{2}\rfloor$.
   Otherwise, $d \leq q \leq \lfloor\frac{k-5}{2}\rfloor$ implies that $N_r(G) \leq \binom{k-q}{r}+(n-k+1)\binom{q}{r-1}= h_r(n, k, q) \leq \max \left\{h_r(n, k, d),h_r\left(n, k, \left\lfloor\frac{k-5}{2}\right\rfloor\right)\right\}$, a contradiction.

\noindent\textbf{ (i)}  $k$ is odd.

\noindent\textbf{Case 1.} $q=\frac{k-1}{2}$.

   Let $T_{1}$ be the set of vertices in $V\left(G^{\prime}\right)$ with degree at least $\frac{k+1}{2}$, i.e.,
$$
T_{1}=\left\{u \in V\left(G^{\prime}\right): d_{G^{\prime}}(u) \geq \frac{k+1}{2}\right\} .
$$
 Then $T_1$ is a clique in $G^{\prime}$.  Let $T_1'=V\left(G^{\prime}\right) \backslash T_{1}$.
   By Lemma \ref{bipartite}, $G^{\prime}\left[T_{1}, T_1'\right]$ is a complete bipartite graph.
We will show that $\left|T_{1}\right|=\frac{k-1}{2}$ or $\left|T_{1}\right|=\frac{k-3}{2}$.

\noindent\textbf{Claim 1.} $\left|T_{1}\right|\leq\frac{k-1}{2}$

\noindent\textbf{Proof.}
Otherwise, $\left|T_{1}\right|\geq\frac{k+1}{2}$.
Since $G^{\prime}\left[T_1, T_1^{\prime}\right]$ is a complete bipartite graph, all vertices in $T'$ with degree at least  $\frac{k+1}{2}$.
   It implies that $T'_1$ is an empty set.
   Thus $G'$ is a  complete graph.
   Since $n\geq k+1$, $G'$ contains a linear forest of size $k$, a contradiction. \qed

\noindent\textbf{Claim 2.} $\left|T_{1}\right|\geq\frac{k-3}{2}$.

 \noindent\textbf{Proof.}  Otherwise, $\left|T_{1}\right|\leq\frac{k-5}{2}$.
   Suppose $|T_{1}|=\frac{k-1}{2}-t$, then $2\leq t\leq \frac{k-1}{2}$.
   Since $G^{\prime}\left[T_1, T_1^{\prime}\right]$ is a complete bipartite graph, the maximum degree of $G^{\prime}\left[T_{1}^{\prime}\right]$ is at most $t$.
   Moreover, $G^{\prime}\left[T_{1}^{\prime}\right]$ is $\mathcal{L}_{2t+1}$-free.
   Otherwise we will find a linear forest of size at least $k$ in $G^{\prime}$.
   By Lemma \ref{k2}, $e(T^{\prime})\leq g(2t,t)\leq 2t(t-1)$ when $t\geq 3$ and $e(T^{\prime})\leq g(2t,t)\leq 6$ when $t=2$.
   Suppose $uv\in E(G'[T_1'])$.
   Since the degrees of $u$ and $v$ are at most $\frac{k-1}{2}$,  $u$ and $v$ have at most $\frac{k-3}{2}$ common neighbors. 	   Thus the edge $uv$ is contained in at most $\binom{\frac{k-3}{2}}{r-2}$ $r$-cliques.

   If $t=2$, then
  $$
\begin{aligned}
  N_r\left(G^{\prime}\right)
  & = N_r(T_1) + N_r^1\left(T_1, T_1'\right) + N_r^{\geq 2}\left(T_1,T_1'\right)\\[1mm]
  & \leq \binom{\frac{k-5}{2}}{r} + \left(n-\frac{k-5}{2}\right)\binom{\frac{k-5}{2}}{r-1}+  6\binom{\frac{k-3}{2}}{r-2}\\[1mm]
  & = \binom{\frac{k-5}{2}}{r} + \left(n-\frac{k+5}{2}\right)\binom{\frac{k-5}{2}}{r-1}+5\binom{\frac{k-5}{2}}{r-1}+6 \binom{\frac{k-3}{2}}{r-2}\\[1mm]
  %& < \binom{\frac{k-5}{2}}{r} + \left(n-\frac{k+5}{2}\right)\binom{\frac{k-5}{2}}{r-1}+5\binom{\frac{k-5}{2}}{r-1}+6 %\binom{\frac{k-3}{2}}{r-2}+4\binom{\frac{k-5}{2}}{r-2}+3\binom{\frac{k-3}{2}}{r-3}+\binom{\frac{k-1}{2}}{r-3}\\
  & < \binom{\frac{k+5}{2}}{r} + \left(n-\frac{k+5}{2}\right)\binom{\frac{k-5}{2}}{r-1}\\[1mm]
  & = h_r\left(n, k, \left\lfloor\frac{k-5}{2}\right\rfloor\right),
\end{aligned}
$$ where the last inequality follows from (\ref{3.1}), a contradiction.

   If $3\leq t\leq \frac{k-1}{2}$, then
   $$
\begin{aligned}
  N_r\left(G^{\prime}\right)
  & = N_r(T_1) + N_r^1\left(T_1, T_1'\right) + N_r^{\geq 2}\left(T_1,T_1'\right)\\[1mm]
  & \leq \binom{\frac{k-1}{2}-t}{r} + \left(n-\frac{k-1}{2}+t\right)\binom{\frac{k-1}{2}-t}{r-1}+ 2t(t-1) \binom{\frac{k-3}{2}}{r-2}\\[1mm]
  & \leq \binom{\frac{k-7}{2}}{r} + \left(n-\frac{k-7}{2}\right)\binom{\frac{k-7}{2}}{r-1}+ \frac{(k-1)(k-3)}{2} \binom{\frac{k-3}{2}}{r-2}\\[1mm]
  %& = \binom{\frac{k-7}{2}}{r} + \left(n-\frac{k+5}{2}+6\right)\left(\binom{\frac{k-5}{2}}{r-1}-\binom{\frac{k-7}{2}}{r-2}\right)+ \frac{(k-1)(k-3)}{2} \binom{\frac{k-1}{2}}{r-2}\\
  & = \binom{\frac{k-7}{2}}{r} + \left(n-\frac{k+5}{2}\right)\left(\binom{\frac{k-5}{2}}{r-1}-\binom{\frac{k-7}{2}}{r-2}\right) +6\binom{\frac{k-7}{2}}{r-1} + \frac{(k-1)(k-3)}{2} \binom{\frac{k-1}{2}}{r-2}\\[1mm]
  & < \binom{\frac{k+5}{2}}{r} + \left(n-\frac{k+5}{2}\right)\binom{\frac{k-5}{2}}{r-1}\\[1mm]
  & = h_r\left(n, k, \left\lfloor\frac{k-5}{2}\right\rfloor\right),
\end{aligned}
$$
where the third inequality follows from  (\ref{3.1}), $n>k^5$ and $r\leq \lfloor\frac{k-3}{2}\rfloor$, a contradiction. \qed

   By Claim 1 and Claim 2, we have $\left|T_{1}\right|=\frac{k-1}{2}$ or $\left|T_{1}\right|=\frac{k-3}{2}$.
   When $\left|T_{1}\right|=\frac{k-3}{2}$,
since $G^{\prime}\left[T_{1}, T_{1}^{\prime}\right]$ is a complete bipartite graph and all the vertices in $T_{1}^{\prime}$ have degree at most $\frac{k-1}{2}$, it follows that all vertices in $T_{1}^{\prime}$ have degree at most one in $G^{\prime}\left[T_{1}^{\prime}\right]$.
   Therefore, $G'[T_{1}^{\prime}]$ consists of independent edges and isolated vertices.
   We claim there are at most two edges in $G^{\prime}\left[T_{1}^{\prime}\right]$.
   Otherwise, one can find $P_{k-2}\cup 3 P_2$ in $G'$, a contradiction.
   Thus, $G'\subseteq H^+\left(n, k-1, \left\lfloor\frac{k-3}{2}\right\rfloor\right)$.
When $\left|T_{1}\right|=\frac{k-1}{2}$,
since $G^{\prime}\left[T_{1}, T_{1}^{\prime}\right]$ is a complete bipartite graph and vertices in $T_{1}^{\prime}$ have degree at most $\frac{k-1}{2}$, it follows that $T_{1}^{\prime}$ forms an independent set of $G^{\prime}$.
   Then $G^{\prime}$ is isomorphic to $H(n, k, \left\lfloor\frac{k-1}{2}\right\rfloor)$.

\noindent\textbf{Case 2.} $q=\frac{k-3}{2}$.

   Let $T_{2}$ be the set of vertices in $V\left(G^{\prime}\right)$ with degree at least $\frac{k+3}{2}$, i.e.,
$$
T_{2}=\left\{u \in V\left(G^{\prime}\right): d_{G^{\prime}}(u) \geq \frac{k+3}{2}\right\} .
$$
Then $T_2$ is a clique in $G^{\prime}$.   Let $T_{2}^{\prime}=V\left(G^{\prime}\right) \backslash T_{2}$.
   By Lemma \ref{bipartite}, $G^{\prime}\left[T_{2}, T_{2}^{\prime}\right]$ is a complete bipartite graph. We will show that $\left|T_{2}\right|=\frac{k-3}{2}$.

\noindent\textbf{Claim 3.} $\left|T_{2}\right|\leq\frac{k-3}{2}$.

\noindent\textbf{Proof.}
Otherwise, $\left|T_{2}\right|\geq\frac{k-1}{2}$.
   The fact $G^{\prime}\left[T_{2}, T_{2}^{\prime}\right]$ is a complete bipartite graph implies that all vertices in $T_{2}^{\prime}$ have degree at least $\frac{k-1}{2}$.
   Therefore $G'$ has no vertex with degree less than or equal to $\frac{k-3}{2}$, which contradicts to the fact that  $G'$ has $(n-k+\frac{k-3}{2}, \frac{k-3}{2})$-P\'osa property. \qed

\noindent\textbf{Claim 4.} $\left|T_{2}\right|\geq\frac{k-3}{2}$.

\noindent\textbf{Proof.}
   Otherwise, $\left|T_{2}\right|\leq\frac{k-5}{2}$.
Suppose $|T_{2}|= \frac{k-1}{2}-t$, where $2\leq t\leq \frac{k-1}{2}$.
   Since $G^{\prime}\left[T_2, T_2^{\prime}\right]$ is a complete bipartite graph, the maximum degree of $G^{\prime}\left[T_{2}^{\prime}\right]$ is at most $t+1$.
   Moreover, $G^{\prime}\left[T_{2}^{\prime}\right]$ is $\mathcal{L}_{2t+1}$-free.
   Otherwise we will find a linear forest of size at least $k$ in $G^{\prime}$.

   When $t=2$,
   since $G^{\prime}\left[T_{2}, T_{2}^{\prime}\right]$ is a complete bipartite graph, $G^{\prime}\left[T_{2}^{\prime}\right]$ is $\mathcal{L}_5$-free with maximum degree at most 3.
   By Lemma \ref{k2}, $e(T^{\prime})\leq g(4,3)< 10=\binom{5}{2}$.
   Then we have $N_r(G'[T'_{2}])\leq \binom{5}{r}$ from Lemma \ref{cliques}.
   Thus the following inequality holds:
   $$
\begin{aligned}
  N_r\left(G^{\prime}\right)
  & = N_r(T_2) + N_r^1(T_2, T'_2) + \sum\limits_{i=2}^5 N_r^i(T_2, T_2')\\[1mm]
  & \leq \binom{\frac{k-5}{2}}{r} + \left(n-\frac{k-5}{2}\right)\binom{\frac{k-5}{2}}{r-1}+ \sum\limits_{i=2}^5\binom{5}{i} \binom{\frac{k-5}{2}}{r-i}\\[1mm]
%  & \leq \binom{\frac{k-5}{2}}{r} +5\binom{\frac{k-5}{2}}{r-1}+ \sum\limits_{i=2}^5\binom{5}{i} \binom{\frac{k-5}{2}}{r-i}+ \left(n-\frac{k+5}{2}\right)\binom{\frac{k-1}{2}-2}{r-1}\\
  & = \binom{\frac{k+5}{2}}{r} + \left(n-\frac{k+5}{2}\right)\binom{\frac{k-5}{2}}{r-1}\\[1mm]
  & = h_r\left(n, k, \left\lfloor\frac{k-5}{2}\right\rfloor\right),
\end{aligned}
$$
where the second equality follows from  (\ref{3.1}), a contradiction.

   When $3\leq t\leq \frac{k-1}{2}$,
   by Lemma \ref{k2}, $e(T^{\prime})\leq g(2t,t+1)\leq 2t^2$.
   Note each edge in $G'[T_2']$ is contained in at most $\binom{\frac{k-1}{2}}{r-2}$ $r$-cliques.
   Thus we have
   $$
\begin{aligned}
  N_r\left(G^{\prime}\right)
  & = N_r(T_2) + N_r^1\left(T_2, T'_2\right) + N_r^{\geq 2}\left(T_2,T'_2\right)\\[1mm]
  & \leq \binom{\frac{k-1}{2}-t}{r} + \left(n-\frac{k-1}{2}+t\right)\binom{\frac{k-1}{2}-t}{r-1}+ 2t^2 \binom{\frac{k-1}{2}}{r-2}\\[1mm]
  & \leq \binom{\frac{k-7}{2}}{r} + \left(n-\frac{k-7}{2}\right)\binom{\frac{k-7}{2}}{r-1}+ \frac{(k-1)^2}{2} \binom{\frac{k-1}{2}}{r-2}\\[1mm]
  %& = \binom{\frac{k-7}{2}}{r} + \left(n-\frac{k+5}{2}+6\right)\left(\binom{\frac{k-5}{2}}{r-1}-\binom{\frac{k-7}{2}}{r-2}\right)+ \frac{(k-1)^2}{2} \binom{\frac{k-1}{2}}{r-2}\\
  & = \binom{\frac{k-7}{2}}{r} + \left(n-\frac{k+5}{2}\right)\left[\binom{\frac{k-5}{2}}{r-1}-\binom{\frac{k-7}{2}}{r-2}\right] +6\binom{\frac{k-7}{2}}{r-1} + \frac{(k-1)^2}{2} \binom{\frac{k-1}{2}}{r-2}\\[1mm]
  & < \binom{\frac{k+5}{2}}{r} + \left(n-\frac{k+5}{2}\right)\binom{\frac{k-5}{2}}{r-1}\\[1mm]
  & = h_r\left(n, k, \left\lfloor\frac{k-5}{2}\right\rfloor\right),
\end{aligned}
$$
where the third inequality follows from (\ref{3.1}), $n>k^5$ and $r\leq \lfloor\frac{k-3}{2}\rfloor$, a contradiction. \qed

By Claim 3 and Claim 4, we have
   $\left|T_{2}\right|=\frac{k-3}{2}$.
Then $G^{\prime}\left[T_{2}^{\prime}\right]$ must be $\mathcal{L}_3$-free. Otherwise we can find a linear forest of size $k$.
   Moreover, each vertex in $G^{\prime}\left[T_{2}^{\prime}\right]$ has degree at most two.
   Thus $G'[T_{2}^{\prime}]$  is a subgraph of $C_3\cup (n-3)K_1$ or $2P_2\cup (n-4)K_1$.
   It follows that $G^{\prime}$ is a subgraph of $H\left(n, k, \left\lfloor\frac{k-3}{2}\right\rfloor\right)$ or $H^+\left(n, k-1, \left\lfloor\frac{k-3}{2}\right\rfloor\right)$.

   Combining the two cases above, we get that $G$ is a subgraph of $H\left(n, k, \left\lfloor\frac{k-1}{2}\right\rfloor\right)$, $H\left(n, k, \left\lfloor\frac{k-3}{2}\right\rfloor\right)$ or $H^+\left(n, k-1, \left\lfloor\frac{k-3}{2}\right\rfloor\right)$.\qed

\noindent\textbf{ (ii)} $k$ is even.

\noindent\textbf{Case 1.} $q=\frac{k-2}{2}$.

   Let $T_{1}$ be the set of vertices in $V\left(G^{\prime}\right)$ with degree at least $\frac{k+2}{2}$, i.e.,
$$
T_{1}=\left\{u \in V\left(G^{\prime}\right): d_{G^{\prime}}(u) \geq \frac{k+2}{2}\right\} .
$$
 Then $T_1$ is a clique in $G^{\prime}$.  Let $T_1'=V\left(G^{\prime}\right) \backslash T_{1}$.
   By Lemma \ref{bipartite}, $G^{\prime}\left[T_{1}, T_1'\right]$ is a complete bipartite graph. We will show that $\left|T_{1}\right|=\frac{k-2}{2}$ or $\left|T_{1}\right|=\frac{k-4}{2}$.

\noindent\textbf{Claim 5.} $\left|T_{1}\right|\leq\frac{k-2}{2}$.

\noindent\textbf{Proof.}
Otherwise, $\left|T_{1}\right|\geq\frac{k}{2}$.
The fact $G^{\prime}\left[T_{1}, T_{1}^{\prime}\right]$ is a complete bipartite graph implies that all vertices in $T_{1}^{\prime}$ have degree at least $\frac{k}{2}$. Then $G^{\prime}$ has no vertex with degree less than or equal to $\frac{k-2}{2}$, which is a contradiction to the fact that $G^{\prime}$ has $(n-k+\frac{k-2}{2}, \frac{k-2}{2})$-P\'osa property. \qed

\noindent\textbf{Claim 6.} $\left|T_{1}\right|\geq\frac{k-4}{2}$.

\noindent\textbf{Proof.}
Otherwise, $\left|T_{1}\right|\leq\frac{k-6}{2}$. Suppose $|T_{1}|= \frac{k}{2}-t$, then $3\leq t\leq \frac{k}{2}$.
   Since $G^{\prime}\left[T_1, T_1^{\prime}\right]$ is a complete bipartite graph, the maximum degree of $G^{\prime}\left[T_{1}^{\prime}\right]$ is at most $t$.
   Moreover, $G^{\prime}\left[T_{1}^{\prime}\right]$ is $\mathcal{L}_{2t}$-free.
   Otherwise we will find a linear forest of size at least $k$ in $G^{\prime}$.
   By Lemma \ref{k2}, $e(T^{\prime})\leq g(2t-1,t)\leq (2t-1)(t-1)$.

   If $t=3$, then
   $$
\begin{aligned}
  N_r\left(G^{\prime}\right)
  & = N_r(T_1) + N_r^1\left(T_1, T_1'\right) + N_r^{\geq 2}\left(T_1,T_1'\right)\\[1mm]
  & \leq \binom{\frac{k-6}{2}}{r} + \left(n-\frac{k-6}{2}\right)\binom{\frac{k-6}{2}}{r-1}+ 10 \binom{\frac{k-2}{2}}{r-2}\\[1mm]
  & = \binom{\frac{k-6}{2}}{r} + \left(n-\frac{k+6}{2}\right)\binom{\frac{k-6}{2}}{r-1}+6\binom{\frac{k-6}{2}}{r-1}+ 10 \binom{\frac{k-2}{2}}{r-2}\\[1mm]
  & < \binom{\frac{k+6}{2}}{r} + \left(n-\frac{k+6}{2}\right)\binom{\frac{k-6}{2}}{r-1}\\[1mm]
  & = h_r\left(n, k, \left\lfloor\frac{k-5}{2}\right\rfloor\right),
\end{aligned}
$$
where the last inequality follows from (\ref{3.1}).

   If $4\leq t\leq \frac{k}{4}$, then
   $$
\begin{aligned}
  N_r\left(G^{\prime}\right)
  & = N_r(T_1) + N_r^1\left(T_1, T_1'\right) + N_r^{\geq 2}\left(T_1,T_1'\right)\\[1mm]
  & \leq \binom{\frac{k}{2}-t}{r} + \left(n-\frac{k}{2}+t\right)\binom{\frac{k}{2}-t}{r-1}+ (2t-1)(t-1) \binom{\frac{k-2}{2}}{r-2}\\[1mm]
  & \leq \binom{\frac{k-8}{2}}{r} + \left(n-\frac{k-8}{2}\right)\binom{\frac{k-8}{2}}{r-1}+ \frac{(k-1)(k-2)}{2} \binom{\frac{k-2}{2}}{r-2}\\[1mm]
  & = \binom{\frac{k-8}{2}}{r} + \left(n-\frac{k+6}{2}\right)\left(\binom{\frac{k-6}{2}}{r-1}-\binom{\frac{k-8}{2}}{r-2}\right)+7\binom{\frac{k-8}{2}}{r-1}+ \frac{(k-1)(k-2)}{2} \binom{\frac{k-2}{2}}{r-2}\\[1mm]
  & < \binom{\frac{k+6}{2}}{r} + \left(n-\frac{k+6}{2}\right)\binom{\frac{k-6}{2}}{r-1}\\[1mm]
  & = h_r\left(n, k, \left\lfloor\frac{k-5}{2}\right\rfloor\right),
\end{aligned}
$$
where the third inequality holds since (\ref{3.1}),  $n>k^5$ and $r\leq \lfloor\frac{k-3}{2}\rfloor$, a contradiction. \qed

By Claim 5 and Claim 6, we have $\left|T_{1}\right|=\frac{k-4}{2}$ or $\left|T_{1}\right|=\frac{k-2}{2}$.
When $\left|T_{1}\right|=\frac{k-4}{2}$,
since $G^{\prime}\left[T_1, T_1^{\prime}\right]$ is a complete bipartite graph, the maximum degree of $G^{\prime}\left[T_{1}^{\prime}\right]$ is at most two.
   Moreover, $G^{\prime}\left[T_{1}^{\prime}\right]$ is $\mathcal{L}_4$-free.
   Therefore, $G^{\prime}\left[T_{1}^{\prime}\right]$ (without isolated vertices) is a subgraph of $\{C_4,C_3\cup P_2, 3P_2\}$.
   Thus $G$ is a subgraph of $H(n,k,\lfloor\frac{k-3}{2}\rfloor)$, $H^{+}(n,k-1,\lfloor\frac{k-3}{2}\rfloor)$ or $H^{++}(n,k-2,\lfloor\frac{k-3}{2}\rfloor)$.
When $\left|T_{1}\right|=\frac{k}{2}-1$,
since $G^{\prime}\left[T_{1}, T_{1}^{\prime}\right]$ is a complete bipartite graph, $G^{\prime}\left[T_{1}^{\prime}\right]$ is $\mathcal{L}_{2}$-free, i.e. there is at most one edge in $G^{\prime}\left[T_{1}^{\prime}\right]$.
   Thus, $G'\subseteq H\left(n, k, \left\lfloor\frac{k-1}{2}\right\rfloor\right)$.

\noindent\textbf{Case 2.} $q=\frac{k-4}{2}$.

   Let $T_{2}$ be the set of vertices in $V\left(G^{\prime}\right)$ with degree at least $\frac{k+4}{2}$, i.e.,
$$
T_{2}=\left\{u \in V\left(G^{\prime}\right): d_{G^{\prime}}(u) \geq \frac{k+4}{2}\right\} .
$$
 Then $T_2$ is a clique in $G^{\prime}$.  Let $T_{2}^{\prime}=V\left(G^{\prime}\right) \backslash T_{2}$.
   By Lemma \ref{bipartite}, $G^{\prime}\left[T_{2}, T_{2}^{\prime}\right]$ is a complete bipartite graph. We will show that $\left|T_{2}\right|=\frac{k-4}{2}$

\noindent\textbf{Claim 7.} $\left|T_{2}\right|\leq\frac{k-4}{2}$.

\noindent\textbf{Proof.}
Otherwise, $\left|T_{2}\right|\geq\frac{k-2}{2}$.
The fact $G^{\prime}\left[T_{2}, T_{2}^{\prime}\right]$ is a complete bipartite graph implies that all vertices in $T_{2}^{\prime}$ have degree at least $\frac{k-2}{2}$.
   Therefore $G'$ has no vertex with degree less than or equal to $\frac{k-4}{2}$, which contradicts to the fact that $G'$ has $(n-k+\frac{k-4}{2}, \frac{k-4}{2})$-P\'osa property. \qed

\noindent\textbf{Claim 8.} $\left|T_{2}\right|\geq\frac{k-4}{2}$.

\noindent\textbf{Proof.}
Otherwise, $\left|T_{2}\right|\leq\frac{k-6}{2}$.
Suppose $|T_{2}|=\frac{k}{2}-t$, then $3\leq t\leq \frac{k}{2}$.
   Since $G^{\prime}\left[T_2, T_2^{\prime}\right]$ is a complete bipartite graph, the maximum degree of $G^{\prime}\left[T_{2}^{\prime}\right]$ is at most $t+1$.
   Moreover, $G^{\prime}\left[T_{2}^{\prime}\right]$ is $\mathcal{L}_{2t}$-free.
   Otherwise, we will find a linear forest of size at least $k$ in $G^{\prime}$.

 When $t=3$, since $G^{\prime}\left[T_{2}, T_{2}^{\prime}\right]$ is a complete bipartite graph, $G^{\prime}\left[T_{2}^{\prime}\right]$ is $\mathcal{L}_6$-free with maximum degree at most 4.
   By Lemma \ref{k2} (ii), $e(T'_{2})\leq 15=\binom{6}{2}$.
   So $N_r(G'[T'_{2}])\leq \binom{6}{r}$ from Lemma \ref{cliques}.
   Then we have
   $$
\begin{aligned}
  N_r\left(G^{\prime}\right)
  & = N_r(T_2) + N_r^1(T_2, T'_2) + \sum\limits_{i=2}^6 N_r^i(T_2, T'_2)\\[1mm]
  & \leq \binom{\frac{k-6}{2}}{r} + \left(n-\frac{k-6}{2}\right)\binom{\frac{k-6}{2}}{r-1}+ \sum\limits_{i=2}^6\binom{6}{i} \binom{\frac{k}{2}-3}{r-i}\\[1mm]
  %& = \binom{\frac{k-6}{2}}{r} + 6\binom{\frac{k-6}{2}}{r-1}+ \sum\limits_{i=2}^6\binom{6}{i} \binom{\frac{k}{2}-3}{r-i}+\left(n-\frac{k+6}{2}\right)\binom{\frac{k-6}{2}}{r-1}\\
  & = \binom{\frac{k+6}{2}}{r} + \left(n-\frac{k+6}{2}\right)\binom{\frac{k-6}{2}}{r-1}\\[1mm]
  & = h_r\left(n, k, \left\lfloor\frac{k-5}{2}\right\rfloor\right),
\end{aligned}
$$
where the second equality follows from  (\ref{3.1}),   a contradiction.

  When $4\leq t\leq \frac{k}{2}$,
   by Lemma \ref{k2}, $e(T^{\prime})\leq g(2t-1,t+1)\leq (2t-1)t$.
   Thus we have
   $$
\begin{aligned}
  N_r\left(G^{\prime}\right)
  & = N_r(T_2) + N_r^1\left(T_2, T'_2\right) + N_r^{\geq 2}\left(T_2,T'_2\right)\\[1mm]
  & \leq \binom{\frac{k}{2}-t}{r} + \left(n-\frac{k}{2}+t\right)\binom{\frac{k}{2}-t}{r-1}+ (2t-1)t \binom{\frac{k}{2}}{r-2}\\[1mm]
  & \leq \binom{\frac{k-8}{2}}{r} + \left(n-\frac{k-8}{2}\right)\binom{\frac{k-8}{2}}{r-1}+ \frac{k(k-1)}{2} \binom{\frac{k}{2}}{r-2}\\[1mm]
  %& = \binom{\frac{k-8}{2}}{r} + \left(n-\frac{k+6}{2}+7\right)\left(\binom{\frac{k-6}{2}}{r-1}-\binom{\frac{k-8}{2}}{r-2}\right)+ \frac{k(k-1)}{2} \binom{\frac{k}{2}}{r-2}\\
  & = \binom{\frac{k-8}{2}}{r} + \left(n-\frac{k+6}{2}\right)\left[\binom{\frac{k-6}{2}}{r-1}-\binom{\frac{k-8}{2}}{r-2}\right] + 7\binom{\frac{k-8}{2}}{r-1}+\frac{k(k-1)}{2} \binom{\frac{k}{2}}{r-2}\\[1mm]
  & < \binom{\frac{k+6}{2}}{r} + \left(n-\frac{k+6}{2}\right)\binom{\frac{k-6}{2}}{r-1}\\[1mm]
  & = h_r\left(n, k, \left\lfloor\frac{k-5}{2}\right\rfloor\right)
\end{aligned}
$$
where the third inequality follows from (\ref{3.1}),  $n>k^5$ and $r\leq \lfloor\frac{k-3}{2}\rfloor$, a contradiction. \qed

By Claim 7 and Claim 8, we have $\left|T_{2}\right|=\frac{k-4}{2}$.
 Since $G^{\prime}\left[T_{2}, T_{2}^{\prime}\right]$ is a complete bipartite graph, all vertices in $T_{2}^{\prime}$ have degree at least $\frac{k-4}{2}$.
   The $(n-k+\frac{k-4}{2},\frac{k-4}{2})$-p\'osa property  implies that there are at most 4 vertices in $T_{2}^{\prime}$ with degree great than 0.
   Thus $G^{\prime}\left[T_{2}^{\prime}\right]$ is a subgraph of $K_4\cup (n-\frac{k+4}{2})K_1$. Then $G\subseteq H\left(n, k, \left\lfloor\frac{k-3}{2}\right\rfloor\right)$.

    Combining the two cases above, we get that $G$ is a subgraph of $H\left(n, k, \left\lfloor\frac{k-1}{2}\right\rfloor\right)$, $H\left(n, k, \left\lfloor\frac{k-3}{2}\right\rfloor\right)$, $H^+\left(n, k-1, \left\lfloor\frac{k-3}{2}\right\rfloor\right)$  or $H^{++}\big(n,k-2,\lfloor\frac{k-3}{2}\rfloor\big)$. The proof is finished. \qed

%%%%%%%%%%%%%%%%%%%%%%%%%%%%%%%%%%%%%%%%%%%%%%%%%%%%%%%%%%%%%%%%%

\section{The clique version of the stability of Erd\H{o}s-Gallai Theorem}

   Notice that a linear forest with at least $2k + 1$ edges has a matching of size at least $k + 1$. A graph $G$ with $\nu(G) \leq k$ must be $\mathcal{L}_{2k+1}$-free.
  Combining  Theorem \ref{stab2} (i) and further discussions, we obtain Theorem \ref{thm11}.

%\begin{theorem}\label{thm111}
%   Let $G$ be a graph on $n$ vertices with $\delta(G) \geq d$ and $\nu(G) \leq k$.
%   If $n \geq (2k+1)^5$ and
%   $$N_r(G)>\max \left\{h_{r}(n, 2k+1, d),h_{r}(n, 2k+1, k-2)\right\},$$
%   then $G$ is a subgraph of $H(n, 2k+1, k)$ or $H(n, 2k+1, k-1)$.
%\end{theorem}

\noindent\textbf{Proof of Theorem \ref{thm11}.}
   Let $G$ be a graph satisfying the conditions of  Theorem \ref{thm11}.
   Then $G$ is $\mathcal{L}_{2k+1}$-free.
   By Theorem \ref{stab2} (i), if  $G\nsubseteq H^+\left(n, 2k, k-1\right)$, then $G$ is a subgraph of $H(n, 2k+1, k)$ or $H\left(n, 2k+1, k-1\right)$.
  Next we will show that if $G\subseteq H^+\left(n, 2k, k-1\right)$, then $G\subseteq H\left(n, 2k+1, k-1\right)$.

   If $G\subseteq H^+\left(n, 2k, k-1\right)$ and  $G\subseteq H\left(n, 2k+1, k-1\right)$, then we are done.
   Now we suppose that $G\subseteq H^+(n, 2k, k-1)$  and $G\nsubseteq H\left(n, 2k+1, k-1\right)$.

   Note that $H^+(n, 2k, k-1)$ can be viewed as a graph obtained from $H(n, 2k-1, k-1)$
   by adding two independent edges, say $x_1y_1$ and $x_2y_2$.
   If $G\subseteq H^+(n, 2k, k-1)$ but $G\nsubseteq H\left(n, 2k+1, k-1\right)$, then $x_1y_1$ and $x_2y_2$ must be in $E(G)$.
   Let $G_1=G-\{x_1,y_1,x_2,y_2\}$.
   Then $G_1\subseteq H(n-4,2k-1,k-1)$ and
\begin{align}\label{G'}
  N_r(G_1)
  & > h_r(n,2k+1,k-2)-4\binom{k-1}{r-1}-2\binom{k-1}{r-2}\notag\\
  %& = \binom{k+3}{r}+(n-k-3)\binom{k-2}{r-1}-4\binom{k-1}{r-1}-2\binom{k-1}{r-2}\notag\\
  & > \binom{k-1}{r}+(n-k-3)\binom{k-2}{r-1}
\end{align}

   Since $G_1\subseteq H(n-4,2k-1,k-1)$, there exists an independent set $I$ satisfies $|I|=n-k-3$ and $d_{G_1}(v)\leq k-1$ for all $v\in I$.
   Suppose that there are $t$ vertices in $I$ with degree $k-1$.
   Then $t\leq k-2$.
   Otherwise, we can find a $(k-1)$-matching $M$ in $G_1$.
 The $(k-1)$-matching  $M$ together with the edges $x_1y_1$ and $x_2y_2$ form a $(k+1)$-matching in $G$, a contradiction.

\noindent\textbf{Case 1.} $t=0$.

   In this case, all vertices in $I$ have degree at most $k-2$.
   Thus
   $$N_r(G_1)\leq \binom{k-1}{r}+(n-k-3)\binom{k-2}{r-1},$$
   contradicting to  (\ref{G'}).

\noindent\textbf{Case 2.} $1\leq t\leq k-2$.

   There are at most $k-2-t$ vertices in $I$ with degree $k-2$.
   Otherwise, for any $S\subseteq V(G_1)\setminus I$, $|N(S)|\geq |S|$. By Hall's Theorem,
   there exists a $(k-1)$-matching $M$ in $G_1$.   The $(k-1)$-matching  $M$ together with the edges $x_1y_1$ and $x_2y_2$ form a $(k+1)$-matching in $G$, a contradiction.
   Thus
\begin{align*}
   N_r(G_1)
   & \leq \binom{k-1}{r}+t\binom{k-1}{r-1}+(k-2-t)\binom{k-2}{r-1}+(n-k-2)\binom{k-3}{r-1}\\[1mm]
   & < \binom{k-1}{r}+(k-1)\binom{k-1}{r-1}+(n-k-3)\binom{k-3}{r-1}\\[1mm]
   & < \binom{k-1}{r}+(n-k-3)\binom{k-2}{r-1},
\end{align*}
   where the last inequality follows from $n>(2k+1)^5$, which is a contradiction to   (\ref{G'}).

   Thus $G\subseteq H^+(n, 2k, k-1)$ implies $G\subseteq H\left(n, 2k+1, k-1\right)$.
   That is, $G$ is a subgraph of $H(n, 2k+1, k)$ or $H(n, 2k+1, k-1)$, completing the proof.
\qed

%%%%%%%%%%%%%%%%%%%%%%%%%%%%%%%%%%%%%%%%%%%%%%%%%%%%%%%%%%%%%%%%%%

\end{document}